\documentclass[journal]{IEEEtran}
\usepackage{amsmath}
\usepackage{amsfonts}

\newcommand{\cN}{\mbox{$\cal{N}$}}
\newcommand{\cA}{\mbox{$\cal{A}$}}
\newcommand{\cI}{\mbox{$\cal{I}$}}
\newcommand{\cC}{\mbox{$\cal{C}$}}

\newtheorem{lemma}{Lemma}
\newtheorem{corollary}{Corollary}
\newtheorem{theorem}{Theorem}

\title{Perfect single error-correcting codes in the Johnson scheme}
\author{Daniel M. Gordon
\thanks{The author is with the IDA Center for Communications Research,
San Diego, CA 92121-1969 (email: gordon@ccrwest.org)} 
}

\begin{document}

\maketitle

\begin{abstract}
  Delsarte conjectured in 1973 that there are no nontrivial pefect
  codes in the Johnson scheme.  Etzion and Schwartz recently showed
  that perfect codes must be $k$-regular for large $k$, and used this
  to show that there are no perfect codes correcting single errors in
  $J(n,w)$ for $n \leq 50{,}000$.  In this paper we show that there are
  no perfect single error-correcting codes for $n \leq 2^{250}$.
\end{abstract}


\section{Introduction}

The Johnson graph $J(n,w)$ has vertices corresponding to $V^n_w$, the
$w$-subsets of the set $\cN = \{1,2,\ldots,n\}$, with two vertices adjacent if their
intersection has size $w-1$.  

The distance between two $w$-sets is half the size of their symmetric
difference.  The $e$-sphere of a point, the set of all $w$-sets within
distance $e$, has cardinality
$$
\Phi_e(n,w) = \sum_{i=0}^e {w \choose i}{n-w \choose i}.
$$

A code $\cC \subset J(n,w)$ 
is called {\em $e$-perfect}
if the $e$-spheres of all the codewords of $\cC$ form a partition of
$V^n_w$. 
Delsarte \cite{delsarte} conjectured that no nontrivial perfect codes
exist in $J(n,w)$.  


Etzion and Schwartz \cite{es2004} 
introduced the concept of $k$-regular
codes.
In this paper we use their results to improve the lower bound
on the size of a 1-perfect code.
The method of proof will be to look at
the factors of $\Phi_1(w,a)$.
We show that $\Phi_1(w,a)$ is squarefree, and
for each prime $p_i|\Phi_1(w,a)$, there is an integer $\alpha_i$ such that
$p_i^{\alpha_i}$ must be close to $n-w$.  Then we will show that the
$\alpha_i$'s are distinct and pairwise coprime, and the
sum of their reciprocals is close to two.
A computer search for perfect powers in short intervals then shows
that no such codes exist with $n < 2^{250}$.

For the rest of this paper we will deal with the case $e=1$, and 
write $n = 2w+a$.  This may be
done without loss of generality, since the complement of an
$e$-perfect code in $J(n,w)$ is $e$-perfect in $J(n,n-w)$.  
Also, to simplify the statements of theorems, we will assume
throughout the paper that $\cC$ is a nontrivial 1-perfect code in $J(n,w)$.

\section{Regularity of 1-Perfect codes}

In this section we summarize the results of Etzion and Schwartz
\cite{es2004} that
we will need.
Let $\cA$ be a $k$-subset of $\cN = \{1,2,\ldots,n\}$.  
For all $0 \leq i \leq k$, define
$$
{\cal C}_{\cA} (i) = | \{ c \in \cC: |c \cap \cA| = i\}|,
$$
and for each $\cI \subseteq \cA$, define
$$
{\cal C}_{\cA} (\cI) = | \{ c \in \cC: c \cap \cA = \cI\}|.
$$

$\cC$ is {\em $k$-regular} if:
\begin{enumerate}
\item There exist numbers $\alpha(0),\alpha(1),\ldots,\alpha(k)$ such
  that for any $k$-set $\cA$ in $\cN$, ${\cal C}_{\cA} (i) = \alpha(i)$,  for
$i=0,1,\ldots,k$.
\item For any $k$-set $\cA$ in $\cN$, there exist numbers
$\beta_{\cA}(0),\beta_{\cA}(1),\ldots,\beta_{\cA}(k)$ such
  that if $\cI \subseteq \cA$, then 
${\cal C}_{\cA} (\cI) = \beta_{\cA}(|\cI|)$.
\end{enumerate}

Etzion and Schwarz give a necessary condition for a code to be
regular:

\begin{theorem}\label{thm:reg}
If $\cC$
is $k$-regular, then
 \begin{equation}
   \label{eq:regular}
\Phi_1(w,a) = 1+w(w+a) \left| {{2w+a-i}\choose{w+a}}    \right.
 \end{equation}
for $i=0,\ldots,k$.
\end{theorem}

They then show that 1-perfect codes must be highly regular.

\begin{theorem}\label{thm:sigma}
$\cC$
is $k$-regular if the polynomial
\begin{equation}
  \label{eq:sigma}
\sigma_1(w,a,m) = m^2 - (2w+a+1)m + w(w+a)+1  
\end{equation}

has no integer roots for $1 \leq m \leq k$.
\end{theorem}

Let 
$$
L(w,a) = \frac{2w+a+1-\sqrt{(a+1)^2+4(w-1)}}{2}.
$$

The smallest root of (\ref{eq:sigma}) is $L(w,a)$, so

\begin{theorem}\label{thm:normal}
$\cC$ is $k$-regular for any $k<L(w,a)$.
\end{theorem}

This means that we can rule out 1-perfect codes by showing that there
is some $i$ with $0 \leq i  \leq  L(w,a)$ such that
(\ref{eq:regular})
is not satisfied.
$L(w,a)$ is an increasing function of $a$, so

\begin{lemma}\label{lemma:lwa}
  $L(w,a) \geq L(w,0) > w-\lceil \sqrt{w} \rceil$.
\end{lemma}





\begin{lemma}\label{lemma:lwa3}
We have
$$0 < a < w/2.$$ 
\end{lemma}

\begin{proof}
Theorem 13 in \cite{es2004}, which is a strengthening of a theorem
of Roos \cite{roos}, gives
$a < w-3$.  If $a=0$ then $\cC$ is a trivial code.

If $a \geq w/2$, then
$$L(w,a) > L\left(w,\frac{w-7}{2}\right) = w-2,$$
so $\cC$ is $(w-2)$-regular.  
$\cC$ is also $(w-1)$-regular, since
$\sigma_1(w,a,w-1) = a-(w-3) \neq 0$ for $a<w-3$.

Since $\cC$ corrects single errors, any two codewords are at least distance 3
apart in $J(n,w)$.  
Let $\cA$ be a $(w-1)$-set contained in some codeword $c_1$.  
Remove any element of $\cA$ and add one not in $c_1$ to get a new
$(w-1)$-set $\cA'$.  Since $\cC$ is $(w-1)$-regular, there is a codeword
$c_2$ containing $\cA'$, but $c_1$ and $c_2$ have distance 2 in
$J(n,w)$, a contradiction.

\end{proof}

\section{Divisors of $\Phi_1(w,a)$}

We will derive necessary conditions for 1-perfect codes by
looking at possible prime divisors of $\Phi_1(w,a)$.  One
tool will be:

\begin{lemma}\label{lemma:binom} (Kummer)
  Let $p$ be a prime.  The number of times $p$ appears in the
  factorization of $a\choose b$ equals the number of carries when
  adding $b$ to $a-b$ in base $p$.
\end{lemma}

Theorem~\ref{thm:normal} and Lemmas~\ref{lemma:lwa} and
\ref{lemma:binom}  imply 

\begin{corollary}\label{cor:carry}
If $p$ is a prime with $p^k | \Phi_1(w,a)$, then
there are at least $k$ carries when adding $w+a$ to $j=w-i$ for 
$j=
\lceil \sqrt{w} \, \rceil+1 ,\lceil \sqrt{w}\, \rceil+2,\ldots,w$.
\end{corollary}

Let
\begin{equation}
  \label{eq:basep}
  w+a = (r_m,r_{m-1},\ldots,r_1,r_0)_p
\end{equation}
be the base $p$ representation of $w+a$, with $r_m \geq 1$.
Let $l = \lfloor m/2\rfloor$.  

\begin{lemma}\label{lemma:highdigits}
  $r_i = p-1$  for $i=l+1, l+2, \ldots, m$.
\end{lemma}

\begin{proof}
  For any $i$ with $\lceil \sqrt{w}\rceil+1 \leq p^i \leq w$, adding
  $p^i$ to $w+a$ must have a carry by Corollary~\ref{cor:carry}, so
  the lemma follows for $i=l+1,\ldots,m-1$.  To complete the proof, we
  need to show that $w\geq p^m$.  We have
$$w+a \geq p^m + (p-1) p^{m-1} \geq \frac{3}{2} p^m.$$
Since $a<w/2$ by Lemma~\ref{lemma:lwa3}, this implies $w > p^m$.
\end{proof}

\begin{theorem}\label{thm:sqfree}
  $\Phi_1(w,a)$ must be squarefree.
\end{theorem}

\begin{proof}
Adding $p^m$ to $w+a$ has only one carry, so by
Corollary~\ref{cor:carry}
only one power of $p$ divides $\Phi_1(w,a)$. 
\end{proof}


\begin{theorem}\label{thm:short}
  For any prime $p$ dividing $\Phi_1(w,a)$, let
$\alpha = m+1 = \lfloor \log_p (w+a) \rfloor + 1$.
Then
  \begin{equation}
    \label{eq:short}
p^\alpha - \lceil \sqrt{w} \rceil - 1 \leq w+a < p^\alpha
  \end{equation}
\end{theorem}

\begin{proof}
We have $w+a < p^\alpha$ from (\ref{eq:basep}).
By
Lemma~\ref{lemma:highdigits}, we must have $r_i = p-1$ for 
$i = l+1, l+2 ,\ldots, m$.
Let
$$(t_l,t_{l-1},\ldots,t_0)_p$$ be the base $p$ representation of
$\lceil \sqrt{w}\,\rceil$.  
The left inequality of (\ref{eq:short}) is equivalent to
\begin{eqnarray*}
p^\alpha-1-(w+a) & = & 
(p-1-r_l,
\ldots,p-1-r_0)_p  \\
& \leq &
(t_l,t_{l-1},\ldots,t_0)_p = \lceil \sqrt{w} \rceil .
\end{eqnarray*}
If this is not satisfied, let $i$ be the largest integer such that
$p-1-r_i > t_i$. 
The number $(t_l, t_{l-1}, \ldots,t_{i+1},t_i+1,0,\ldots,0)_p$ is
greater than $\lceil \sqrt{w}\rceil$ and has no carries when when
added to $w+a$ in base $p$, which contradicts
Corollary~\ref{cor:carry}.
\end{proof}

Thus we have that $p^\alpha$ is in a short interval around $w+a$.
We will use this result in the following form:

\begin{corollary}\label{cor:bound1}
  For a prime $p$ dividing $\Phi_1(w,a)$, we have
  \begin{equation}
    \label{eq:cor2}
0 <  \log_{w+a} p -\frac{1}{\alpha} <
\frac{1}{\alpha} \left(\frac{1}{\sqrt{w+a}} + 
\frac{4}{(w+a)}\right).
  \end{equation}
\end{corollary}

\begin{proof}
From  (\ref{eq:short}), we have
\begin{eqnarray*}
p^\alpha > w+a & \geq  &
p^{\alpha}\left(1-\frac{\left\lceil\sqrt{w}\right\rceil+1}{p^{\alpha}}\right) \\
& > &  p^\alpha \left( 1-\frac{1}{\sqrt{w+a}} - \frac{2}{w+a} \right)
\end{eqnarray*}
using $\lceil \sqrt{w}\rceil +1 < \sqrt{w+a}+2$.
Taking the log base $w+a$, we have
$$
\alpha \log_{w+a} p >
1 > \alpha \log_{w+a} p + \log_{w+a} \left( 1-\frac{1}{\sqrt{w+a}}
  - \frac{2}{w+a} \right)
$$
Using the bound $ -\log (1-x) < x+x^2$ for $x <1/2$
gives the corollary. 
\end{proof}

\section{Powers in Short Intervals}

Theorem~\ref{thm:short} shows that for a 1-perfect code to exist,
several prime powers must be close to $w+a$.  
Having a large number of prime powers in a short interval seems
unlikely.  Loxton \cite{loxton} showed (a gap in the proof was later
fixed by Bernstein \cite{bernstein}) that the number of perfect powers
in $[w,w+\sqrt{w}\,]$ is at most
$$
\exp (40 \sqrt{\log \log w \log \log \log w}).
$$
Loxton conjectured that the number of perfect
powers in such an interval is bounded by a constant, 
but a proof seems very far off.

For the rest of this paper, take 
\begin{equation}
  \label{eq:primes}
p_1 p_2 \ldots p_r = \Phi_1(w,a) = 1+w(w+a).
\end{equation}


Taking the log of (\ref{eq:primes}) gives
$$
\sum_{i=1}^r \log_{w+a} p_i = \log_{w+a} (w(w+a) + 1),
$$
so
\begin{eqnarray}
  \label{eq:bound2}
0 &<& \sum_{i=1}^r \log_{w+a} p_i - (1 + \log_{w+a} w) \nonumber \\
&= & \log_{w+a}(1+\frac{1}{w(w+a)})\\
&\leq& \frac{1}{w(w+a)}. \nonumber
\end{eqnarray}


\begin{theorem}
\label{thm:egyptian}
$$
\left|\sum_{i=1}^r \frac{1}{\alpha_i}
 - (1 + \log_{w+a} w) \right| <
\frac{4}{\sqrt{w+a}}.
$$

\end{theorem}

\begin{proof}
If $\sum_{i=1}^r \frac{1}{\alpha_i} - (1 + \log_{w+a} w) \geq 0$, then
  the theorem follows immediately from (\ref{eq:bound2}) and
  Corollary~\ref{cor:bound1}.  Otherwise, 
  summing (\ref{eq:cor2}) we have
  \begin{eqnarray*}
0 &< &  (1 + \log_{w+a} w) - \sum_{i=1}^r \frac{1}{\alpha_i} \\
& < & 
\sum_{i=1}^r \log_{w+a} p_i - \sum_{i=1}^r \frac{1}{\alpha_i} \\
& < &
\sum_{i=1}^r \frac{1}{\alpha_i} \left(\frac{1}{\sqrt{w+a}} +
  \frac{4}{w+a}\right)\\
& < & 2\frac{2}{\sqrt{w+a}}.
  \end{eqnarray*}
\end{proof}

Clearly the constant 4 in Theorem~\ref{thm:egyptian} can be strengthened, but this will be enough for
our purposes.

For $0 < a < w/2$, we have $w+a < 3w/2$, so 
$$
1-\log_{w+a} 3/2 < \log_{w+a} w < 1
$$
and Theorem \ref{thm:egyptian} says that we have an Egyptian fraction
representing a number close to 2.  Etzion and Schwartz
showed that there are no 1-perfect codes with $n \leq 50000$, and
so
\begin{equation}
  \label{eq:egypt}
\frac{1}{\alpha_1} + \frac{1}{\alpha_2} + \ldots
\frac{1}{\alpha_r} \in \left[1.934, 2.026\right].
\end{equation}

\begin{lemma}
  The $\alpha_i$'s are distinct and pairwise coprime.
\end{lemma}

\begin{proof}
We cannot have $\alpha_i = \alpha_j = 1$,
since then
$p_i, p_j > (w+a)$ implies $p_i p_j > 1+w(w+a) = \Phi_1(w,a)$, contradicting
(\ref{eq:primes}). 

Suppose we have $\alpha_i$, $\alpha_j$ with $\gcd(\alpha_i,
\alpha_j) = g>1$.
Then by Theorem \ref{thm:short},
$p_i^{\alpha_i}$ and $p_j^{\alpha_j}$ are two $g^{\mbox{{\rm th}}}$
powers 
in an interval around $w+a$ of length $\sqrt{w+a}$, 
which is impossible.

\end{proof}

For an integer $k$, let $p^-(k)$ denote the smallest prime factor of
$k$.  

\begin{corollary}\label{cor:seven}
  Some $\alpha_i$ has $p^-(\alpha_i) \geq 7$.
\end{corollary}
  
\begin{proof}
If there are more than four $\alpha$'s, clearly one of them must have
a prime factor bigger than 5.  For four $\alpha$'s, the set
$\{1,2,3,5\}$ has sum of reciprocals $2.033$, which by
(\ref{eq:egypt}) is too big, and
an easy computation finds that
any set of powers of these numbers has a sum of reciprocals that is too
small.  The largest is $\{1,2,3,25\}$, with sum $1.8733$.
\end{proof}

Let $\gamma(n)$ denote the largest squarefree divisor of $n$.
The $abc$ conjecture asserts that, for any $\epsilon>0$ 
there are only finitely many integers $a$, $b$ and $c$ such
that $a+b=c$ and 
$$
\max\{a,b,c\} \leq C_\epsilon \gamma(abc)^{1+\epsilon}.
$$
See \cite{abc} for information and references about the $abc$ conjecture

For any choice of $\alpha$'s satisfying (\ref{eq:egypt}), 
Masser-Oesterl\'{e}'s $abc$ conjecture implies there are
only a finite number of solutions.  
For example, take 
$\alpha_1 = 1$, 
$\alpha_2 = 2$, 
$\alpha_3 = 3$, and 
$\alpha_4 = 7$.
Let $a=p_3^3$, $c=p_4^7$ and $b$ be their difference, which is at most
$\max \{p_3^{3/2},p_4^{7/2}\}$
by Theorem~\ref{thm:short}.  Then
\begin{eqnarray*}
\max\{a,b,c\}  \approx w+a &\leq & C_\epsilon p_3 p_4 c \\
&<& (w+a)^{(1+\epsilon)(1/3+1/7+1/2)} \\
&<& (w+a)^{0.98}
\end{eqnarray*}
for all but finitely many $w$'s.





\section{A New Lower Bound for $n$}

While we cannot show that there are no perfect codes,
Theorem~\ref{thm:short} gives us an efficient way to search for possible
codes, by searching for powers in short intervals.  

To show a bound of $2^C$ for $n$, we need to check for 
primes $a,b \geq 2$ and integers $3 \leq p,q<C$ with
$$
0 < a^p-b^q < \sqrt{a^p}.  
$$
It suffices to consider prime values of $p$ and $q$, since any $k$th
power is also a $p^-(k)$th power.
It is possible to run
through the possibilities efficiently.  Let $\{p_1=3,p_2 =5,\ldots,p_k\}$
be the odd primes up to $C$.
The following procedure
will find all pairs $i$, $j$ and integers $b_i$, $b_j$ for which 
$b_i^{p_i}$ and $b_j^{p_j}$ are
close:

\begin{enumerate}
\item Start with $b_1 = b_2 \cdots = b_k = 2$.  Compute powers
$c_i = b_i^{p_i}$ for $i=1,2,\ldots,k$.
\item Let $c_i$ be the smallest power, and $c_j$ the second smallest.  
 Compare them to see if they are close enough.
\item Increment the base $b_i$, recompute $c_i$, and continue.
\item Stop when all powers are larger than $2^C$.
\end{enumerate}

If two powers less than $2^C$ are in a short interval, they will
eventually be the two smallest powers in the list, and will be found.
A heap (see, for example, \cite{knuth3}) 
is an efficient data structure to maintain
the powers in, requiring only one comparison to find the two smallest
powers, and $\leq \log_2 k$ steps to reorder the heap after changing
$c_i$.

Note that the above algorithm looks for any integers $b_i$ and $b_j$
with powers in a short interval, not just primes.  Only considering primes
would reduce the number of comparisons, but complicate the
rule for stepping the bases $b_i$.

In five hours on a 2.6 GHz Opteron, 
an implementation of this algorithm
eliminated everything up to $2^{109}$.  It found 60 powers higher than
squares in short intervals, most of which involved a cube and fifth
power.  By Corollary~\ref{cor:seven}, we may discount these.
The only higher powers are
given in Table~\ref{tab:pow}.  

\begin{table}
\begin{center}
\begin{tabular}{|c|c|c|}
\hline
$p_1^{\alpha_1}$ & $p_2^{\alpha_2}$ & difference
\\\hline
$2^7$ & $5^3$ & $3$ \\
$13^3$ & $3^7$ & $10$ \\ 
$3251^3$ & $32^7$ & $83883$ \\ 
$33^7$ & $3493^3$ & $178820$ \\ 
$1965781^3$ & $498^7$ & $1539250669$ \\ 
\hline\end{tabular}
\end{center}
\caption{Pairs of Higher Powers in Short Intervals up to $2^{109}$}\label{tab:pow}
\end{table}

Only the first two pairs are powers of primes, and they are in the
range already
eliminated by Etzion and Schwartz's result.  The larger ones all
involve at least one composite, so they do not result in a 1-perfect
code.  Therefore we have

\begin{theorem}
  There are no 1-perfect codes in $J(n,w)$ for all $n<2^{109}$.
\end{theorem}

Finally, we may bootstrap this result to a stronger one.  Using this
larger bound in Theorem \ref{thm:egyptian}, we can tighten (\ref{eq:egypt})
to
$$
\frac{1}{\alpha_1} + \frac{1}{\alpha_2} + \ldots
\frac{1}{\alpha_r} \in \left[1.99, 2.001\right].
$$

No set of four $\alpha_i$'s have a sum of reciprocals in this
interval, and the only sets of five that do are 
$\{1,2,3,7,k\}$, where $k\in [41,71]$ with $\gcd(k,2\cdot 3 \cdot 7)=1$.
Any set of six $\alpha_i$'s clearly have two $\alpha$'s with a factor
$\geq 7$, so we have

\begin{corollary}
  At least two $\alpha_i$'s have $p^-(\alpha_i) \geq 7$.
\end{corollary}


Therefore we may do a search as above, but
starting with $p_1=7$ instead of 3.  
The search work is proportional to $2^{C/p_1}$,
so this greatly reduces the search time.
A search for seventh and higher powers up to $2^{250}$ in a short
interval took four
hours and 
found none, so 

\begin{theorem}
  There are no 1-perfect codes in $J(n,w)$ for all $n<2^{250}$.
\end{theorem}

\vspace{.2in}

\noindent
{\bf Acknowledgments.}
The author would like to thank the anonymous referee, who suggested
changes which greatly improved the presentation of this paper, and
pointed out Lemma~\ref{lemma:lwa3}.

\bibliographystyle{plain}

  

\end{document}